\documentclass[12pt,a4paper]{article}
\usepackage[cp850]{inputenc}
\usepackage{latexsym,amsfonts}

\newtheorem{Thm}{Theorem}[section]

\newtheorem{Con}[Thm]{Conjecture}

\newtheorem{tab}[Thm]{Table}
\newtheorem{note}[Thm]{Note}

\begin{document}

\title{Maximal $(120,8)$-arcs in projective planes
of order 16 and related designs} 

\author{  Vladimir D. Tonchev\thanks{ Research supported by NSA Grant
H98230-16-1-0011.}  and Tim Wagner\thanks{ Research supported by NSA Grant
H98230-16-1-0011.} \\ 
 Michigan Technological University,
           Houghton, MI 49931, USA\\
           tonchev@mtu.edu, tjwagner@mtu.edu }

\maketitle

\begin{abstract}

The resolutions and maximal sets of compatible resolutions 
of all 2-$(120,8,1)$ designs arising from  
maximal $(120,8)$-arcs 
in the known projective planes of order 16 are computed. 
It is shown that each of these designs is embeddable  
in a unique way in a projective plane of order 16.



\end{abstract}

\section{Introduction}

A 2-$(v,k,\lambda)$ design (or shortly, a 2-design) is a pair $D$=$\{ X, \cal{B} \}$
of a set $X$ of $v$ {\it points} and a collection of subsets
of $X$ of size $k$ called {\it blocks}, such that every two points appear together in
exactly $\lambda$ blocks \cite{BJL}, \cite{CRC}.
Every point of a
$2$-$(v,k,\lambda)$ design is contained in
$r=\lambda(v-1)/(k-1)$ blocks, and the total number
of blocks is $b=v(v-1)\lambda/k(k-1)$.

The incidence matrix of a design $D$ is a $(0,1)$-matrix $A=(a_{ij})$ with rows
labeled by the blocks, columns labeled by the points, where
$a_{i,j}=1$ if the $i$th block contains the $j$th point, and $a_{i,j}=0$
otherwise. If $p$ is a prime, the $p$-rank of a design $D$ is the rank
of its incidence matrix of the finite field of order $p$.

Two designs are {\it isomorphic} if there is a bijection between their
point sets that maps every block of the first design to a block of the second
design. An {\it automorphism} of a design is any isomorphism of the design
to itself.
The set of all automorphisms of $D$ form the automorphism group
$Aut(D)$ of $D$.

The {\it dual} design $D^\perp$ of a design $D$
has as points the blocks of $D$, and as blocks the points of $D$.
A 2-$(v,k,\lambda)$ design is {\it symmetric} if $b=v$,
or equivalently,  $r=k$.
The dual design $D^\perp$ of a symmetric 2-$(v,k,\lambda)$ design $D$
is a symmetric design with the same parameters as $D$.
A symmetric design $D$ is {\it self-dual} if $D$ and $D^{\perp}$ are
isomorphic.

A design with $\lambda=1$ is called a {\it Steiner} design.
An affine plane of order $n$ ($n\ge 2$),
 is a Steiner 2-$(n^2, n, 1)$ design.
A projective plane of order $n$ is a symmetric
Steiner 2-$(n^2 + n +1, n+1, 1)$ design with $n\ge 2$.
The classical (or Desarguesian) plane $PG(2,p^t)$
of order $n=p^t$, where $p$ is prime and $t\ge 1$, has as points
the 1-dimensional subspaces of the 3-dimensional vector space $V_3$
over the finite field of order $p^t$, and as blocks (or {\it lines}), the
2-dimensional subspaces of $V_3$.

Let $D=\{ X, \cal{B} \}$ be a Steiner 2-$(v,k,1)$ design with
point set $X$, collection of blocks $\cal{B}$, and let $v$ be a multiple
of $k$,  $v=nk$. Since every point of $X$ is contained in
$r=(v-1)/(k-1)=(nk-1)/(k-1)$ blocks, it follows that $k-1$ divides $n-1$.
Thus, $n-1=s(k-1)$ for some integer $s\ge 1$, and
\[ v=nk = (sk -s +1)k. \]
A {\it parallel class} $P$ is a set of $v/k=n$
pairwise disjoint blocks,
and a {\it resolution} of $D$ is a partition of the collection of blocks
$\cal{B}$ into $r=(v-1)/(k-1)=sk+1$ parallel classes.
A design is {\it resolvable} if it admits a resolution.

Any 2-$((sk -s +1)k, k, 1)$ design with $s=1$ is 
an affine plane of order $k$, and admits exactly one resolution.
If $s>1$, a resolvable 2-$((sk -s +1)k, k, 1)$ design may admit more than one
resolution.

Following \cite{ton-res}, we call two resolutions
 $R_1$, $R_2$,
\begin{equation}
\label{def}
 R_1 = P^{(1)}_1 \cup P^{(1)}_2 \cup \cdots P^{(1)}_r, \  R_2 =
 P^{(2)}_1 \cup P^{(2)}_2 \cup \cdots P^{(2)}_r
\end{equation}
{\em compatible} if they share one parallel class,
$ P^{(1)}_i = P^{(2)}_j$, 
and
\[ |P^{(1)}_{i'} \cap P^{(1)}_{j'}|\le 1 \]
for  $i' \neq i$ and  $j' \neq j$.

More generally, a set of $m$ resolutions
$ R_{1}, \ldots, R_{m}$
 is  compatible
if every two of these resolutions are  compatible.

Suppose that $\cal{P}$ is a projective plane of order $q=sk$. A {\it maximal}
$\{ (sk-s+1)k, k \}$-{\it arc} \cite{Hir},
is a set $\cal{A}$ of $(sk-s+1)k$ points
of $\cal{P}$ such that every line  of $\cal{P}$ is ether disjoint from $\cal{A}$
or meets $\cal{A}$ in exactly $k$ points.
The collection of lines of $\cal{P}$ which have no points in common with $\cal{A}$
determines a maximal  $\{ (sk-k+1)s, s)$-arc ${\cal{A}}^{\perp}$ in the dual plane
${\cal{P}}^{\perp}$.

Maximal arcs with $1<k < q$ do not exist in any Desarguesian plane of odd order
$q$ \cite{BBM},
and are known to exist in every Desarguesian plane of
order $q=2^t$ \cite{Den}, for any $k=2^i, \ 1\le i<t$,
 as well as in some non-Desarguesian planes of even order
\cite{HST}, \cite{PRS}.

If $k>1$, the non-empty intersections of
a maximal $\{ (sk-s+1)k, k \}$-arc $\cal{A}$ with the lines of
a projective plane $\cal{P}$ of order $q=sk$
are the blocks of a resolvable 2-$((sk-s+1)k,k,1)$ design $D$.
Similarly, if $s>1$, the corresponding $\{ (sk-k+1)s, s \}$-arc ${\cal{A}}^{\perp}$
in the dual plane is the point set of a resolvable
2-$((sk-k+1)s, s, 1)$ design $D^{\perp}$.
We will refer to $D$ (resp. $D^{\perp}$) as a design embeddable in
$\cal{P}$ (resp. ${\cal{P}}^\perp$) as a maximal arc.
The points of $D^{\perp}$ determine a set of $(sk-k+1)s$ mutually compatible
resolutions of $D$. Respectively, the points of $D$ determine a set
of $(sk-s+1)k$ mutually compatible resolutions of $D^{\perp}$.

Two maximal arcs ${\cal{A}'}$, ${\cal{A}''}$ in a projective plane
$\cal{P}$ are {\it equivalent }
if there is an automorphism of $\cal{P}$ that maps
  ${\cal{A}'}$ to  ${\cal{A}''}$.
We note that the designs associated with equivalent arcs
are necessarily isomorphic, while the converse is not true
in general.

In \cite{ton-res}, one of the authors of this paper
 proved the following upper bound
on the number of pairwise compatible resolutions of a 
2-$((sk-s+1)k, k, 1)$ design.

\begin{Thm}
\label{th}
Let $S=\{ R_1,\ldots, R_m \}$ be a set of $m$ mutually  compatible resolutions
of a  2-$((sk -s +1)k, k, 1)$ design $D=\{ X, \cal{B} \}$.
Then
\[ m\le (sk-k+1)s. \]
The equality
\[ m=(sk-k+1)s \]
holds if and only if there exists a projective plane $\cal{P}$ of order $sk$
such that $D$ is embeddable in $\cal{P}$ as a maximal $\{ (sk-s+1)k, k \}$-arc.

\end{Thm}

This paper summarizes the computation of all parallel classes, resolutions,
and compatible sets of resolutions of maximum size
of the 2-$(120,8,1)$ associated with maximal (120,8)-arcs 
in projective planes of order 16.
The main result can be formulated as follows.

\begin{Thm}
\label{main}
 Every 2-$(120,8,1)$ design associated with a maximal
$(120,8)$-arc in a known projective plane $\cal{P}$ of order 16 admits exactly one
compatible set of resolutions of maximal size, meeting the bound of Theorem \ref{th},
and consequently, is uniquely embeddable in $\cal{P}$.\\
\end{Thm}

\section{Maximal (120,8)-arcs and related 2-$(120,8,1)$ designs}

There are 22 nonisomorphic projective planes of order 16 that are known currently.
Four planes, $PG(2,16)$, SEMI2, SEMI4, and BBH1 (in the notation of \cite{PRS})
are self-dual, and there are nine planes which are not self-dual:
HALL, LMRH, JOWK, DSFP, DEMP, BBH2, JOHN, BBS4, and MATH \cite{PRS}.
Lists with the collections of lines of these planes are available
at Eric Moorhouse's web page \cite{M}.

In \cite{PRS}, Penttila, Royle, and Simpson enumerated and classified 
up to equivalence all hyperovals in the known planes of order 16.
A hyperoval $\cal{A}$
 in a plane $\cal{P}$ of order 16 is a maximal (18,2)-arc, and
its dual arc ${\cal{A}}^{\perp}$ is a maximal (120,8)-arc in the dual
plane ${\cal{P}}^{\perp}$.
Since two maximal arcs ${\cal{A}'}$ ${\cal{A}''}$ in a plane $\cal{P}$
are equivalent if and only if their dual arcs $({\cal{A}'})^{\perp}$,
$({\cal{A}''})^{\perp}$ are equivalent, the results from \cite{PRS} imply
the classification of all maximal (120,8)-arcs in the known projective
planes of order 16, up to equivalence.

We used the data about the inequivalent hyperovals
 graciously provided to the authors by Gordon F. Royle,
to compute the corresponding dual (120,8)-arcs and the related 2-$(120,8,1)$
designs.
 The 93 inequivalent hyperovals give rise to 93 inequivalent (120,8)-arcs.
For each 2-$(120,8,1)$ design $D$ associated with an arc in the dual plane
of the plane containing the corresponding hyperoval, we computed
all parallel classes of $D$, all resolutions of $D$, and all
compatible sets of maximal size 18. The parallel classes were found as
13-cliques in a graph $\Gamma$ having as vertices the blocks of $D$,
where two blocks are adjacent in $\Gamma$ if they are disjoint.
The resolutions were computed as 17-cliques 
in a graph $\Delta$ having as vertices
the parallel classes of $D$, where two parallel classes are adjacent
in $\Delta$ if they do not share any block. Finally, we computed the compatible
sets of maximal size 18 as 18-cliques in a graph $\cal{E}$ 
having as vertices the resolutions of $D$,
where adjacency is defined according to the definition of compatible resolutions
given in the preceding section.
For these computations, we wrote algorithms using  Magma
\cite{Magma} and Cliquer \cite{cliquer}.

The results of these computations are summarized in Table \ref{tab1}.

\begin{note}
{\rm
The number of parallel classes ranges from 153 to 221, 
while the number of resolutions is 18 in all but one
notable exception, namely, the 2-$(120,8,1)$ design corresponding to the dual
(120,8)-arc of the regular hyperoval in $PG(2,16)$, with  group
of order 16,320. The number of resolutions of this particular design is 137 
(cf. Table \ref{tab1}). However,
among those 137 resolutions, there is only one set of 18 pairwise compatible
resolutions, thus, this design, as well as all remaining  designs,
are embeddable in a unique way in a projective plane of order 16.
}
\end{note}

\begin{note}
{\rm
The 2-ranks of the 2-$(120,8,1)$ designs associated
with maximal (120,8)-arcs, as well as their groups,
were previously computed by Laurel Carpenter \cite{C-th}.
The 2-ranks of the 93 designs range from 65 to 94,
and the minimum 65 is achieved
only by the two designs in the Desarguesian plane, $PG(2,16)$,
corresponding to the regular hyperoval with a
group of order 16,320, and the
Lunelli-Sce hyperoval \cite{LS}
(also known as the Hall hyperoval \cite{Hall}),
with a group of order 144
(cf. Table \ref{tab2}).
This supports a conjecture from \cite{C-th}, stating that
the 2-rank of any design
associated with a hyperoval in $PG(2,2^t)$ is $3^t - 2^t$,
as well as the following stronger conjecture formulated in \cite{ton-res},
which generalizes a conjecture by A. E. Brouwer \cite{Br}.

}
\end{note}

\begin{Con}
\label{conj}
If $\cal{D}$ is a 2-$(2^{2t-1}-2^{t-1}, 2^{t-1},1)$ design ($t \ge 2$),
with an incidence matrix $A$, then
\[ rank_{2}(A) \ge 3^t - 2^t, \]
and the equality
\[   rank_{2}(A) = 3^t - 2^t \]
holds if and only if $\cal{D}$ is embeddable as a maximal $(2^{2t-1}-2^{t-1}, 2^{t-1})$-arc
in $PG(2^t, 2)$.
\end{Con}

Conjecture \ref{conj} is trivially true for $t=2$, and its validity for $t=3$
follows from the results of \cite{MTW}.

\begin{center}
\begin{tab}
\label{tab1}
\end{tab}
\begin{tabular}{|c|c|c|c|c|c|}
\hline {\bf Hyperoval} & {\bf Plane} &  {\bf $|$Aut(\emph{D})$|$} & {\bf 2-rank} & {\bf \# Par. Cl.} & {\bf \# Resolutions} 
\\\hline 1 & $PG(2,\,16)$& 144 & 65 & 153  &  18 
\\\hline 2 & $PG(2,\,16)$& 16320 & 65 & 221  & 137 
\\\hline 1 & SEMI2 & 3 & 81 & 153  &  18 
\\\hline 2 & SEMI2 & 8 & 81 & 153  &  18 
\\\hline 3 & SEMI2 & 8 & 81 & 153  &  18 
\\\hline 4 & SEMI2 & 16 & 81 & 153  &  18 
\\\hline 5 & SEMI2 & 16 & 81 & 153  &  18 
\\\hline 6 & SEMI2 & 16 & 80 & 153  &  18 
\\\hline 7 & SEMI2 & 16 & 81 & 153  &  18 
\\\hline 8 & SEMI2 & 16 & 80 & 153  &  18 
\\\hline 9 & SEMI2 & 16 & 81 & 153  &  18 
\\\hline 10 & SEMI2 & 16 & 80 & 153  &  18 
\\\hline 11 & SEMI2 & 16 & 81 & 153  &  18 
\\\hline
\end{tabular}
\end{center}
\pagebreak
\begin{center}
\begin{tabular}{|c|c|c|c|c|c|}
\hline {\bf Hyperoval} & {\bf Plane} &  {\bf $|$Aut(\emph{D})$|$} & {\bf 2-rank} & {\bf \# Par. Cl.} & {\bf \# Resolutions}
\\\hline 12 & SEMI2 & 16 & 81 & 153  &  18 
\\\hline 13 & SEMI2 & 16 & 80 & 153  &  18 
\\\hline 14 & SEMI2 & 16 & 81 & 153  &  18 
\\\hline 15 & SEMI2 & 16 & 81 & 153  &  18 
\\\hline 16 & SEMI2 & 32 & 81 & 157  &  18 
\\\hline 17 & SEMI2 & 32 & 80 & 157  &  18 
\\\hline 1 & SEMI4 & 16 & 81 & 153  &  18 
\\\hline 2 & SEMI4 & 16 & 81 & 153  &  18 
\\\hline 3 & SEMI4 & 16 & 81 & 153  &  18 
\\\hline 1 & HALL & 16 & 80 & 153  &  18 
\\\hline 2 & HALL & 64 & 81 & 157  &  18 
\\\hline 3 & HALL & 64 & 80 & 157  &  18 
\\\hline 4 & HALL & 320 & 80 & 173  &  18 
\\\hline 1 & HALL.d & 2 & 81 & 153  &  18 
\\\hline 2 & HALL.d & 2 & 81 & 153  &  18 
\\\hline 3 & HALL.d & 6 & 81 & 153  &  18 
\\\hline 1 & LMRH & 16 & 83 & 153  &  18 
\\\hline 2 & LMRH & 16 & 86 & 153  &  18 
\\\hline 3 & LMRH & 16 & 86 & 153  &  18 
\\\hline 4 & LMRH & 16 & 86 & 153  &  18 
\\\hline 5 & LMRH & 64 & 86 & 157  &  18 
\\\hline 6 & LMRH & 112 & 82 & 153  &  18 
\\\hline 1 & LMRH.d & 14 & 89 & 153  &  18 
\\\hline 1 & JOWK & 16 & 82 & 153  &  18 
\\\hline 2 & JOWK & 16 & 82 & 153  &  18 
\\\hline 3 & JOWK & 16 & 83 & 153  &  18 
\\\hline 4 & JOWK & 16 & 82 & 153  &  18 
\\\hline 5 & JOWK & 64 & 82 & 157  &  18 
\\\hline 6 & JOWK & 112 & 82 & 153  &  18 
\\\hline 1 & JOWK.d & 14 & 83 & 153  &  18 
\\\hline 1 & DSFP & 16 & 86 & 153  &  18 
\\\hline 2 & DSFP & 16 & 86 & 153  &  18 
\\\hline 3 & DSFP & 16 & 86 & 153  &  18 
\\\hline 4 & DSFP & 16 & 86 & 153  &  18 
\\\hline 5 & DSFP & 16 & 86 & 153  &  18 
\\\hline 6 & DSFP & 16 & 86 & 153  &  18 
\\\hline
\end{tabular}
\end{center}
\pagebreak
\begin{center}
\begin{tabular}{|c|c|c|c|c|c|}
\hline {\bf Hyperoval} & {\bf Plane} &  {\bf $|$Aut(\emph{D})$|$} & {\bf 2-rank} & {\bf \# Par. Cl.} & {\bf \# Resolutions}
\\\hline 7 & DSFP & 16 & 86 & 153  &  18 
\\\hline 8 & DSFP & 16 & 86 & 153  &  18 
\\\hline 9 & DSFP & 16 & 86 & 153  &  18 
\\\hline 10 & DSFP & 16 & 86 & 153  &  18 
\\\hline 11 & DSFP & 16 & 86 & 153  &  18 
\\\hline 12 & DSFP & 16 & 86 & 153  &  18 
\\\hline 13 & DSFP & 16 & 86 & 153  &  18 
\\\hline 14 & DSFP & 16 & 86 & 157  &  18 
\\\hline 15 & DSFP & 16 & 86 & 153  &  18 
\\\hline 16 & DSFP & 16 & 85 & 153  &  18 
\\\hline 17 & DSFP & 16 & 86 & 153  &  18 
\\\hline 18 & DSFP & 16 & 86 & 153  &  18 
\\\hline 19 & DSFP & 16 & 86 & 153  &  18
\\\hline 20 & DSFP & 16 & 86 & 153  &  18 
\\\hline 21 & DSFP & 16 & 86 & 153  &  18 
\\\hline 22 & DSFP & 64 & 86 & 157  &  18 
\\\hline 1 & DEMP & 16 & 85 & 153  &  18 
\\\hline 2 & DEMP & 16 & 84 & 153  &  18 
\\\hline 3 & DEMP & 16 & 85 & 153  &  18 
\\\hline 4 & DEMP & 16 & 85 & 153  &  18 
\\\hline 5 & DEMP & 16 & 83 & 153  &  18 
\\\hline 6 & DEMP & 16 & 85 & 153  &  18 
\\\hline 7 & DEMP & 16 & 84 & 153  &  18 
\\\hline 8 & DEMP & 16 & 85 & 153  &  18 
\\\hline 9 & DEMP & 16 & 84 & 153  &  18 
\\\hline 10 & DEMP & 16 & 85 & 153  &  18 
\\\hline 11 & DEMP & 16 & 84 & 153  &  18 
\\\hline 12 & DEMP & 64 & 83 & 157  &  18 
\\\hline 13 & DEMP & 64 & 83 & 157  &  18 
\\\hline 1 & DEMP.d & 2 & 85 & 153  &  18 
\\\hline 2 & DEMP.d & 6 & 85 & 153  &  18 
\\\hline 1 & JOHN & 16 & 94 & 157  &  18 
\\\hline 1 & BBS4 & 16 & 94 & 157  &  18 
\\\hline 1 & BBH2 & 4 & 94 & 153  &  18 
\\\hline 2 & BBH2 & 4 & 94 & 153  &  18 
\\\hline 1 & MATH & 8 & 90 & 153  &  18 
\\\hline
\end{tabular}
\end{center}
\pagebreak
\begin{center}
\begin{tabular}{|c|c|c|c|c|c|}
\hline {\bf Hyperoval} & {\bf Plane} &  {\bf $|$Aut(\emph{D})$|$} & {\bf 2-rank} & {\bf \# Par. Cl.} & {\bf \# Resolutions}
\\\hline 1 & MATH.d & 4 & 91 & 153  &  18 
\\\hline 2 & MATH.d & 4 & 91 & 153  &  18 
\\\hline 3 & MATH.d & 8 & 92 & 153  &  18 
\\\hline 1 & BBH1 & 8 & 90 & 153  &  18 
\\\hline 2 & BBH1 & 16 & 92 & 153  &  18 
\\\hline 3 & BBH1 & 32 & 90 & 157  &  18 
\\\hline
\end{tabular}
\end{center}

\vspace{5mm}

\begin{center}
\begin{tab}{2-ranks of 2-$(120,8,1)$ designs}
\label{tab2}
\end{tab}
\end{center}
\begin{tabular}{|c|r|r|r|r|r|r|r|r|r|r|r|r|r|}
\hline
2-rank & \ 65 \ & \ 80 \ & \ 81 \ & \ 82 \ & \ 83 \ & \ 84 \ & \ 85 \ & \ 86 \
 & \ 89 \ & \ 90 \ & \ 91 \ & \ 92 \ & \ 94 \ \\
\hline
Frequency & \ 2 \ & \ 8 \ & \ 19 \  & \ 6 \ & \ 7 \ & \ 4 \ & \ 10 \ & \ 25 \
 & \ 1 \ & \ 3 \ & \ 2 \ & \ 2 \ & \ 4 \ \\
\hline
\end{tabular}

\section{Acknowledgments}

The authors wish to thank Gordon F. Royle for graciously providing
the files with the data concerning the hyperovals in the known
projective planes of order 16,
and Tim Penttila for the helpful discussion and for bringing
to our attention some relevant references.

Research supported by NSA Grant  H98230-16-1-0011.


\begin{thebibliography}{100}

\bibitem{BBM} S. Ball, A. Blokhuis, and F. Mazzocca, Maximal arcs in Desarguesian
planes of odd order do not exist, {\it Combinatorica}, {\bf 17} (1997), 31-41.

 
\bibitem{BJL} T. Beth, D. Jungnickel,  H. Lenz, \emph{Design Theory}, 2nd edition.
Cambridge University Press, Cambridge (1999).

\bibitem{Magma} W.  Bosma, J. Cannon, {\it Handbook of Magma Functions},
Department of Mathematics, University of Sydney, 1994.

\bibitem{Br}  A. E. Brouwer, Some unitals and their embeddings
in projective planes of order 9, in: {\it Geometries and Groups},
M. Aigner and D. Jungnickel, eds., {\it Lecture Notes in Mathematics},
{\bf 893} (1981), 183-188.

\bibitem{C96} L. L. Carpenter, Oval designs in Desarguesian projective planes,
{\it Designs, Codes, and Cryptography}, {\bf 9} (1996), 51-59.

\bibitem{C-th} L. L. Carpenter, {\it Designs and codes from hyperovals},
PHD Dissertation, Clemson University, 1996.

\bibitem{CRC} C. J. Colbourn, J. F. Dinitz, eds., {\it Handbook
of Combinatorial Designs}, Second Edition, Chapman \& Hall/CRC,
Boca Raton, 2007.

\bibitem{Den} R. H. F. Denniston, Some maximal arcs in finite projective
planes, {\it J. Combin. Theory}, {\bf 6} (1969), 317-319.

\bibitem{Hall} M. Hall, Jr.,  Ovals in the Desarguesian plane of order 16,
{\it Ann. Mat. Pura Appl.} (4) {\bf 102} (1975), 159–176.

\bibitem{HST} N. Hamilton, S. D. Stoichev, and V. D. Tonchev,
Maximal arcs and disjoint maximal arcs in projective planes of order 16,
{\it J. Geometry} {\bf 67} (2000), 117-126.

\bibitem{Hir} J.~W.~P. Hirschfeld,
\emph{Projective Geometries over Finite Fields (2nd edition)},
Oxford University Press (1998).

\bibitem{LS}  L. Lunelli, M. Sce,
K-archi completi nei piani proiettivi desarguesiani di rango 8 e 16,
Centro di Calcoli Numerici, Politecnico di Milano, Milan 1958,
pp. 15.

\bibitem{MTW}  G. McGuire, V. D. Tonchev, H. N. Ward,
Characterizing the Hermitian and Ree unitals on 28 points,
{\it Designs, Codes and Cryptography}, {\bf 13} (1998), 57-61.

\bibitem{M} E. Moorhouse, Projective Planes of Order 16,\\
http://www.uwyo.edu/moorhouse/pub/planes16/.


\bibitem{cliquer}
S. Niskanen, P. R. J. \"Osterg\r ard, Cliquer User's Guide, Version 1.0. Tech. Rep. T48,
Communications Laboratory, Helsinki University of Technology, Espoo, Finland, 2003.

\bibitem{PRS} T. Penttila, G. F. Royle, and M. K. Simpson,
Hyperovals in the known projective planes of order 16,
{\it J. Combin. Designs}, {\bf 4}, No. 1 (1996), 59-65.


\bibitem{ton-res} V. D. Tonchev,  On resolvable Steiner 2-designs and maximal arcs
  in projective planes,
Designs, Codes, and Cryptography,
DOI: 10.1007/s10623-016-0304-6,
to appear,  http://link.springer.com/article/10.1007/s10623-016-0304-6,
http://arxiv.org/abs/1607.06785.
 
\end{thebibliography}
\end{document}